\input amstex
\input amsppt.sty
\magnification=\magstep1
\hsize=30truecc
\vsize=22.2truecm
\baselineskip=16truept
\TagsOnRight
\pageno=1
\nologo
\def\Z{\Bbb Z}
\def\N{\Bbb N}

\def\l{\left}
\def\r{\right}
\def\bg{\bigg}
\def\({\bg(}
\def\[{\bg\lfloor}
\def\){\bg)}
\def\]{\bg\rfloor}
\def\t{\text}
\def\f{\frac}

\def\p{\ (\roman{mod}\ p)}

\def\bi{\binom}
\def\eq{\equiv}

\def\ls{\leqslant}
\def\gs{\geqslant}
\def\mo{\roman{mod}}
\def\ord{\roman{ord}}

\def\ve{\varepsilon}

\def\da{\delta}

\def\Proof{\noindent{\it Proof}}

\def\Remark{\medskip\noindent{\it  Remark}}

\hbox {Int. J. Number Theory 7(2011), no.\,3, 645--662.}
\bigskip
\topmatter
\title On some new congruences for binomial coefficients\endtitle
\author Zhi-Wei Sun$^1$ and Roberto Tauraso$^2$\endauthor
\leftheadtext{Zhi-Wei Sun and Roberto Tauraso}
\rightheadtext{On some new congruences for binomial coefficients}
\affil $^1$Department of Mathematics, Nanjing University\\
 Nanjing 210093, People's Republic of China\\zwsun\@nju.edu.cn
\\{\tt http://math.nju.edu.cn/$\sim$zwsun}
\medskip
$^2$Dipartimento di Matematica
\\Universit\`a di Roma ``Tor Vergata"
\\Roma 00133, Italy
\\tauraso\@mat.uniroma2.it
\\{\tt http://www.mat.uniroma2.it/$\sim$tauraso}
\endaffil
\abstract In this paper we establish some new congruences involving
central binomial coefficients as well as Catalan numbers. Let $p$ be
a prime and let $a$ be any positive integer. We determine
$\sum_{k=0}^{p^a-1}\bi{2k}{k+d}$ mod $p^2$ for $d=0,\ldots,p^a$ and
$\sum_{k=0}^{p^a-1}\bi{2k}{k+\delta}$ mod $p^3$ for $\delta=0,1$. We
also show that
$$\f1{C_n}\sum_{k=0}^{p^a-1}C_{p^an+k}\eq1-3(n+1)\l(\f{p^a-1}3\r)\ (\mo\ p^2)$$
for every $n=0,1,2,\ldots,$
where $C_m$ is the Catalan number $\bi{2m}m/(m+1)$, and $(\f{\cdot}3)$
is the Legendre symbol.
\endabstract
\thanks 2010 {\it Mathematics Subject Classification}.\,Primary 11B65;
Secondary 05A10,\,11A07.
\newline\indent The initial draft of this paper was written in
2007 and posted as {\tt arXiv:0709.1665}.
\newline\indent The first author is the corresponding author, and he
is supported by the National Natural Science Foundation (grant
10871087) and the Overseas Cooperation Fund (grant 10928101) of
China.
\endthanks
\endtopmatter
\document

\heading{1. Introduction}\endheading

For $n\in\N=\{0,1,2,\ldots\}$, the $n$th Catalan number is given by
$$C_n=\f{1}{n+1}\bi{2n}n=\bi{2n}n-\bi{2n}{n+1}.$$
Here is an alternate definition:
$$C_0=1\ \ \t{and}\ \ C_{n+1}=\sum_{k=0}^nC_kC_{n-k}\ \ (n=0,1,2,\ldots).$$
The Catalan numbers play important roles in combinatorics; they
arise naturally in many enumeration problems (see, e.g., [St, pp.\,219--229]).
For example, $C_n$ is the number of
binary parenthesizations of a string of $n+1$ letters,
and it is also the number of ways to triangulate a convex $(n+2)$-gon into $n$ triangles
by $n-1$ diagonals that do not intersect in their interiors.

In 2006 H. Pan and Z. W. Sun [PS] employed a useful identity to
deduce many congruences on Catalan numbers, in particular they
determined the partial sum $\sum_{k=0}^{p-1}C_k$ modulo a prime $p$
in terms of the Legendre symbol $(\f{\cdot}3)$. For any $a\in\Z$,
$(\f a3)\in\{0,\pm1\}$ satisfies the congruence $a\eq(\f a3)\ (\mo\ 3)$.

In this paper we establish some further congruences involving
Catalan numbers and related central binomial coefficients.

For an assertion $A$, we adopt Iverson's notation:
$$[A]=\cases1&\t{if}\ A\ \t{holds},\\0&\t{otherwise}.\endcases$$
For $d\in\N$ we set
$$S_d=\sum_{0<k<d}\f{(-1)^{k-1}}k\l(\f{d-k}3\r).\tag1.1$$
Note that
$$S_0=S_1=0,\ S_2=1,\ S_3=-\f 32,\ S_4=\f 56,\ S_5=\f 5{12},\ S_6=-\f{21}{20}.$$

 Here is our first theorem.

\proclaim{Theorem 1.1} Let $p$ be a prime, and let
$d\in\{0,\ldots,p^a\}$ with $a\in\Z^+=\{1,2,3,\ldots\}$. Set
$$\aligned F(d)=&\f12\(\sum_{k=0}^{p^a-1}\bi{2k}{k+d}-\l(\f{p^a-d}3\r)\)
\\&+(-1)^{p}p^aS_d-p[p=3]\l(\f d3\r).\endaligned\tag1.2$$
Then
$$F(d)\eq p[p=2\ \&\ 3\nmid (-1)^a+d]\ (\mo\ p^2)\tag1.3$$
and
$$F(p^a-d)\eq F(d)\ (\mo\ p^3).\tag1.4$$
\endproclaim

\Remark\ 1.1. Let $p$ be any prime, and let $a\in\Z^+$ and
$\da\in\{0,1\}$. (1.4) in the case $d=\da$ yields the congruence
$$\aligned&\f12\(\sum_{k=0}^{p^a-1}\bi{2k}{k+\da}-\l(\f{p^a-\da}3\r)\)
\\\eq&2\da p[p=3]+(-1)^{p-1}p^a\sum_{k=1}^{p^a-1}\f{(-1)^k}k\l(\f{p^a-\da-k}3\r)\ (\mo\ p^3).
\endaligned\tag1.5$$

\proclaim{Corollary 1.1} Let $p$ be any prime and let $a\in\Z^+$. For $d\in\{0,1,\ldots,p^a\}$, we have
$$\sum_{k=0}^{p^a-1}\bi{2k}{k+d}\eq\l(\f {p^a-d}3\r)-p[p=3]\l(\f d3\r)+2p^aS_d\ (\mo\ p^2).\tag1.6$$
Also,
$$\sum_{k=0}^{p^a-1}C_k\eq1-3\l(\f{p^a-1}3\r)\eq\f{3(\f{p^a}3)-1}2\ (\mo\ p^2)\tag1.7$$
and
$$\sum_{k=0}^{p^a-1}kC_k\eq\l(\f{p^a-1}3\r)-p[p=3]\eq\f{1-(\f{p^a}3)}2\ (\mo\ p^2).\tag1.8$$
\endproclaim
\Proof. (1.6) holds since $2F(d)\eq0\ (\mo\ p^2)$. As
$C_k=\bi{2k}k-\bi{2k}{k+1}$ and $kC_k=\bi{2k}{k+1}$, both (1.7) and (1.8)
follow from (1.6) with $d=0,1$. \qed

\medskip
\Remark\ 1.2. Let $p$ be a prime and let $a\in\Z^+$. For
$d=0,1,\ldots,p^a$, (1.6) implies that
$$\sum_{k=0}^{p^a-1}\bi{2k}{k+d}\eq\l(\f {p^a-d}3\r)\ (\mo\ p),$$
which was proved by Pan and Sun [PS, Theorem 1.2] in the case $a=1$
via a sophisticated combinatorial identity. (1.6) in the case
$d=0,1$ yields that
$$\sum_{k=0}^{p^a-1}\bi{2k}k\eq\l(\f{p^a}3\r)\ (\mo\ p^2)\tag1.9$$
and $$ \sum_{k=1}^{p^a-1}\bi{2k}{k+1}\eq\l(\f{p^a-1}3\r)-p[p=3]\
(\mo\ p^2).\tag1.10$$ (1.9) in the case $a=1$ implies the following
observation of A. Adamchuk [A] (who told the second author that he
could not find a proof): If $p>3$ then
 $$\sum_{k=1}^{p+(\f{p+1}3)}\bi{2k}k\eq0\ (\mo\ p^2).$$
 (Recall the Wolstenholme congruence
 $\f12\bi{2p}p=\bi{2p-1}{p-1}\eq1\ (\mo\ p^3)$ for $p>3$ (see, e.g., [HT]).)
Recently, (1.9) in the case $a=2$  was posed by D. Callan as a
problem in [C]; in fact, (1.9) in the case $a\in\{2,3,4\}$ was also
observed by A. Adamchuk [A] slightly earlier who could not provide a
proof.
\medskip

Now we state our second theorem.
\proclaim{Theorem 1.2} Let $p$ be a prime. Let $d\in\{0,1,\ldots,p^a\}$ with $a\in\Z^+$, and let
$m,n\in\N$ with $m\gs n$. Then
$$\aligned&\f{n+1}{\bi mn}\sum_{k=0}^{p^a-1}\bi{p^am+2k}{p^an+k+d}-(n+1)\l(\f{p^a-d}3\r)-(m-n)\l(\f d3\r)
\\\eq&((m-n)^2+(m+1)(n+2))p^aS_d-[p=3]\l(\f d3\r)p(n+1)(m+n+1)
\\&+[p=2\ \&\ 3\nmid d-(-1)^a]pm(n+1)\ (\mo\ p^2).
\endaligned\tag1.11$$
In particular,
$$\aligned&\f1{C_n}\sum_{k=0}^{p^a-1}\bi{2p^an+2k}{p^an+k+d}
-n\l(\f d3\r)-(n+1)\l(\f{p^a-d}3\r)
\\\eq&p^a(n+1)(3n+2)S_d-[p=3]p(n+1)\l(\f d3\r)\ (\mo\ p^2).
\endaligned\tag1.12$$
\endproclaim

\proclaim{Corollary 1.2} Let $p$ be a prime and let $a,n\in\N$ with $a>0$. Then
$$\f1{C_n}\sum_{k=0}^{p^a-1}C_{p^an+k}\eq 1-3(n+1)\l(\f{p^a-1}3\r)\ (\mo\ p^2)\tag1.13$$
and
$$\aligned&\f1{C_n}\sum_{k=0}^{p^a-1}kC_{p^an+k}+[p=3]p(n+1)
\\\eq&(1-p^a)n+(3p^an+1)(n+1)\l(\f{p^a-1}3\r)\ (\mo\ p^2).
\endaligned\tag1.14$$
\endproclaim

\Remark\ 1.3. Note that (1.13) and (1.14) are extensions of (1.7)
and (1.8).

\medskip

(1.7) and (1.9) in the case $a=1$ suggest the following open
problem.

\proclaim{Problem 1.1} Are there
any composite numbers $n\not\eq0\ (\mo\ 3)$ such that
$$\sum_{k=0}^{n-1}\bi{2k}k\eq\l(\f n3\r)\ (\mo\ n^2)\ ?$$
Are there any composite numbers $n\not\eq0\ (\mo\ 3)$ satisfying
$$\sum_{k=0}^{n-1}C_k\eq1-3\l(\f{n-1}3\r)\ (\mo\ n^2)\ ?$$
\endproclaim

\Remark\ 1.4. It seems that the answers to Problem 1.1 are negative.
We have confirmed this for $n\ls 10^4$ via Maple.

\medskip

 We are going to do some preparations  in the next section.
 We will show Theorem 1.1 in Sections 3, and  prove Theorem 1.2 and Corollary 1.2 in Section 4.

\heading{2. Some lemmas}\endheading

As usual, for a prime $p$ and an integer $m$, we define
$$\ord_p(m)=\sup\{a\in\N:\,p^a\mid m\}$$
(and thus $\ord_p(0)=+\infty$). A congruence modulo $+\infty$
refers to the corresponding equality.

\proclaim{Lemma 2.1} Let $p$ be any prime, and let $a\in\Z^+$ and $m,n\in\N$
with $m\gs n$. Then
$$\bi{p^am}{p^an}\bigg/\bi mn\eq1+[p=2]pn(m-n)\ (\mo\ p^{2+\ord_p(n)}).\tag2.1$$
\endproclaim
\Proof. (2.1) holds trivially when $n$ is $0$ or $m$.

Below we assume $0<n<m$. Observe that
$$\bi{p^am}{p^an}=\prod_{j=0}^{p^an-1}\f{p^am-j}{p^an-j}=\prod_{i=0}^{n-1}\f{p^am-p^ai}{p^an-p^ai}\times
\prod\Sb 0\ls j<p^an\\p^a\nmid j\endSb\l(1+\f{p^a(m-n)}{p^an-j}\r).$$
Thus
$$\align\f{\bi{p^am}{p^an}}{\bi mn}\eq&1+(m-n)\sum\Sb 0\ls j<p^an\\p^a\nmid j\endSb\f{p^a}{p^an-j}
=1+(m-n)\sum\Sb 0<i<p^an \\p^a\nmid i\endSb\f{p^a}i
\\\eq&1+(m-n)\sum_{q=0}^{n-1}\sum_{k=1}^{p^a-1}\f{p^a}{p^aq+k}
\eq1+(m-n)n\sum_{k=1}^{p^a-1}\f{p^a}k
\\\eq&1+[p=2]p(m-n)n\ (\mo\ p^{2+\ord_p(m-n)})
\endalign$$
since
$$2\sum_{k=1}^{p^a-1}\f{p^a}k=\sum_{k=1}^{p^a-1}\l(\f{p^a}k+\f{p^a}{p^a-k}\r)
=\sum_{k=1}^{p^a-1}\f{p^a}k\cdot\f{p^a}{p^a-k}\eq0\ (\mo\ p^2)$$
and
$$\sum_{k=1}^{2^a-1}\f{2^a}k-\sum_{j=1}^{2^{a-1}-1}\f{2^{a-1}}j=\sum^{2^a-1}\Sb k=1\\2\nmid k\endSb\f{2^a}k\eq0\ (\mo\ 2^a).$$
Similarly,
$$\align&\bi{pm}{pn}\bigg/\bi mn=\bi{pm}{p(m-n)}\bigg/\bi m{m-n}
\\\eq&1+[p=2]pn(m-n)\ (\mo\  p^{2+\ord_p(n)}).
\endalign$$
We are done. \qed

\Remark\ 2.1. By a deep result of Jacobsthal (see, e.g., [Gr]),
if $p>3$ is a prime and $m\gs n\gs0$ are integers, then
$\bi{pm}{pn}=\bi mn\l(1+p^3mn(m-n)v\r)$ for some $p$-adic integer $v$.
However, the proof of Jacobsthal's result is complicated and not easily found in modern literature.
\medskip

\proclaim{Lemma 2.2} Let $p$ be a prime and let $a\in\Z^+$. Then
$$\aligned&\f12\sum_{k=1}^{p^a-1}\bi{2k}k+\sum_{k=1}^{p^a-1}\bi{2k}{k+1}
\\=&\bi{2p^a-1}{p^a-1}-1\eq p[p=2]+p^2[p=3]\ (\mo\ p^3).\endaligned\tag2.2$$
\endproclaim
\Proof. Clearly
$$\bi{2k}k+\bi{2k}{k+1}=\bi{2k+1}{k+1}=\f12\bi{2k+2}{k+1}$$
for all $k\in\N$. Thus
$$2\sum_{k=0}^{p^a-1}\bi{2k}k+2\sum_{k=0}^{p^a-1}\bi{2k}{k+1}
=\sum_{k=0}^{p^a-1}\bi{2k+2}{k+1}=\sum_{k=1}^{p^a}\bi{2k}k$$ and
hence
$$\f12\sum_{k=1}^{p^a-1}\bi{2k}k+\sum_{k=1}^{p^a-1}\bi{2k}{k+1}=\f12\bi{2p^a}{p^a}-1=\bi{2p^a-1}{p^a-1}-1.$$

By Lemma 2.1,
$$\f{\bi{2p^i}{p^i}}{\bi{2p^{i-1}}{p^{i-1}}}\eq1+[p=2]p^{2i-1}\eq1\
(\mo\ p^{i+1})\quad\t{for}\ i=2,3,\ldots.$$
So we have
$$\f12\bi{2p^a}{p^a}\eq\f12\bi{2p}p=\bi{2p-1}{p-1}\eq1+p[p=2]+p^2[p=3]\ (\mo\ p^3)$$
by applying the well-known Wolstenholme congruence in the last step.

Clearly (2.2) follows from the above.
\qed

\proclaim{Lemma 2.3} Let $n>1$ and $d$ be integers. Then
$$\sum_{k=1}^{n-1}\bi nk\l(\f{d+k}3\r)=(1+(-1)^n-3[3\mid n])\l(\f{d-n}3\r).\tag2.3$$
\endproclaim
\Proof. Let $\omega$ denote the cubic root $(-1+\sqrt{-3})/2$ of unity.
As observed by E. Lehmer [L1] in 1938, for any $r\in\Z$ we have
$$\align 3\sum_{k\eq r\,(\mo\ 3)}\bi nk=&\sum_{k=0}^n\bi nk(1+\omega^{k-r}+\omega^{2(k-r)})
\\=&2^n+\omega^{-r}(1+\omega)^n+\omega^{-2r}(1+\omega^2)^n
\\=&2^n+\omega^{-r}(-\omega^2)^n+\omega^r(-\omega)^n
\\=&2^n+(-1)^n(\omega^{n+r}+\omega^{-n-r})
\\=&\cases2^n+2(-1)^n&\t{if}\ 3\mid n+r,\\2^n-(-1)^n&\t{if}\ 3\nmid n+r.
\endcases
\endalign$$
It follows that
$$\align & \sum_{k=1}^{n-1}\bi{n}k\l(\f{d+k}3\r)+\l(\f d3\r)+\l(\f{d+n}3\r)
\\=&\sum_{k=0}^n\bi nk\l(\f{d+k}3\r)=\sum_{k\eq 1-d\,(\mo\ 3)}\bi{n}k-\sum_{k\eq 2-d\,(\mo\ 3)}\bi{n}k
\\=&\f{2^{n}+(-1)^n(3[3\mid n+1-d]-1)}3-\f{2^{n}+(-1)^n(3[3\mid n+2-d]-1)}3
\\=&(-1)^n([3\mid n+1-d]-[3\mid n+2-d])=(-1)^n\l(\f{d-n}3\r).
\endalign$$
Since
$$\l(\f{d-n}3\r)+\l(\f d3\r)+\l(\f{d+n}3\r)=3[3\mid n]\l(\f d3\r),$$
we finally obtain the desired (2.3). \qed

\Remark\ 2.2. The evaluation of $\sum_{k\eq r\,(\mo\ 12)}\bi nk$
with $n\in\N$ and $r\in\Z$, can be found in [Su].

\proclaim{Lemma 2.4} Let $p$ be a prime and let $a\in\Z^+$. Then
$$p^a\(3\sum^{p^a-1}\Sb k=1\\3\mid k-p^a\endSb\f{(-1)^k}k-\sum_{k=1}^{p^a-1}\f{(-1)^k}k\)\eq p[p=2]+\f p2[p=3]\ (\mo\ p^3).\tag2.4$$
\endproclaim
\Proof. By Lehmer's result mentioned in the proof of Lemma 2.3,
$$\align&3\sum\Sb 0<k<p^a\\3\mid k-p^a\endSb\bi{p^a}k+3[3\mid p^a]+3
=3\sum_{k\eq p^a\,(\mo\ 3)}\bi{p^a}k
\\=&2^{p^a}+(-1)^{p^a}(3[3\mid p^a+p^a]-1)=2^{p^a}-3[3\mid p^a]-(-1)^{p^a}.
\endalign$$
Therefore
$$\align&\sum_{k=1}^{p^a-1}\bi{p^a}k-3\sum^{p^a-1}\Sb k=1\\3\mid k-p^a\endSb\bi{p^a}k
\\=&(1+1)^{p^a}-2-\l(2^{p^a}-6[p=3]-3-(-1)^p\r)=6[p=3]+2[p=2].
\endalign$$

For $k=1,\ldots,p^a-1$, since $p^a/k\eq0\ (\mo\ p)$ and
$$\bi{p^a-1}{k-1}(-1)^{k-1}=\prod_{0<j<k}\l(1-\f{p^a}j\r)\eq1-\sum_{0<j<k}\f{p^a}j\ (\mo\ p^2),$$
we have
$$\bi{p^a}k+p^a\f{(-1)^k}k=\f{p^a}k\(\bi{p^a-1}{k-1}+(-1)^{k}\)\eq(-1)^k\f{p^a}k\sum_{0<j<k}\f{p^a}j\ (\mo\ p^3).$$
If $1\ls k\ls p^a-1$ and $p^{a-1}\nmid k$, then $p^a/k\eq0\ (\mo\
p^2)$. Note also that $-1\eq 1\ (\mo\ 2)$. Thus
$$\align&\sum_{k=1}^{p^a-1}(3[3\mid k-p^a]-1)(-1)^k\f{p^a}k\sum_{0<j<k}\f{p^a}j
\\\eq&\sum_{k=1}^{p-1}(3[3\mid p^{a-1}k-p^a]-1)(-1)^{p^{a-1}k}\f{p^a}{p^{a-1}k}\sum_{0<j<p^{a-1}k}\f{p^a}j
\\\eq&\sum_{k=1}^{p-1}(3[3\mid p^{a-1}(k-p)]-1)(-1)^k\f pk\sum_{0<j<k}\f{p^a}{p^{a-1}j}
\\\eq&\sum_{k=1}^{p-1}(3[p\ne3\ \&\ 3\mid k-p]-1)(-1)^k\f pk\sum_{0<j<k}\f pj\ (\mo\ p^3).
\endalign$$

Combining the above,
$$p^a\sum_{k=1}^{p^a-1}(3[3\mid k-p^a]-1)\f{(-1)^k}k\ \mo\ p^3$$
does not depend on $a$. So it suffices to prove (2.4) for $a=1$.

When $a=1$ and $p\in\{2,3\}$, (2.4) can be verified directly.

Below we assume $p>3$. With the help of
Wolstenholme's result $\sum_{k=1}^{p-1}1/k\eq0\ (\mo\ p^2)$ (cf. [Gr] or [HT]), we see that
$$\align &3\sum^{p-1}\Sb k=1\\3\mid k-p\endSb\f{(-1)^k}k-\sum_{k=1}^{p-1}\f{(-1)^k}k
\\\eq&-\sum_{j=1}^{(p-1)/2}\f1j+3\sum_{j=1}^{\lfloor p/3\rfloor}\f1{p-3j}
-6\sum_{j=1}^{\lfloor p/6\rfloor}\f1{p-6j}\ (\mo\ p^2).
\endalign$$
Recall the following congruences of E. Lehmer [L2]:
$$\align \sum_{j=1}^{(p-1)/2}\f1j\eq&-2q_p(2)+pq_p^2(2)\ (\mo\ p^2),
\\\sum_{j=1}^{\lfloor p/3\rfloor}\f1{p-3j}\eq&\f{q_p(3)}2-\f p4q_p^2(3)\ (\mo\ p^2),
\\\sum_{j=1}^{\lfloor p/6\rfloor}\f1{p-6j}\eq&\f{q_p(2)}3+\f{q_p(3)}4
-\f p6q_p^2(2)-\f p8q_p^2(3)\ (\mo\ p^2),
\endalign$$
where $q_p(2)=(2^{p-1}-1)/p$ and $q_p(3)=(3^{p-1}-1)/p$ are Fermat quotients.
Consequently,
$$\sum_{k=1}^{p-1}(3[3\mid k-p]-1)\f{(-1)^k}k\eq0\ (\mo\ p^2)$$
and hence (2.4) holds for $a=1$.

In view of the above, we have completed the proof. \qed

\heading{3. Proof of Theorem 1.1}\endheading

In this section we prove Theorem 1.1 on the basis of Section 2.

\proclaim{Lemma 3.1} Let $p$ be a prime and let $d\in\{1,\ldots,p^a-1\}$
with $a\in\Z^+$. For the function $F$ given by $(1.1)$ we have
$$\aligned &F(d-1)+F(d)+F(d+1)
\\\eq&2\bi{p^a}d-p^a\(\f{(-1)^{p^a-d}}d+(-1)^{p^a-1}\f{(-1)^{d}}{p^a-d}\)\ (\mo\ p^3)
\\\eq&0\ (\mo\ p^2).\endaligned\tag3.1$$
\endproclaim
\Proof. Observe that
$$\align&\bi{2k}{k+d-1}+\bi{2k}{k+d}+\bi{2k}{k+d+1}
\\=&\bi{2k+1}{k+d}+\(\bi{2k+1}{k+d+1}-\bi{2k}{k+d}\)
=\bi{2k+2}{k+d+1}-\bi{2k}{k+d}\endalign$$
for every $k=0,1,2,\ldots$. Thus
$$\align&\sum_{k=0}^{p^a-1}\(\bi{2k}{k+d-1}+\bi{2k}{k+d}+\bi{2k}{k+d+1}\)
\\=&\sum_{k=0}^{p^a-1}\(\bi{2(k+1)}{(k+1)+d}-\bi{2k}{k+d}\)=\bi{2p^a}{p^a+d}-\bi{0}{d}=\bi{2p^a}{p^a-d}.
\endalign$$
Note also that
$$\align&\sum_{\ve=-1}^1\sum_{0<k<d+\ve}\f{(-1)^k}k\l(\f{d+\ve-k}3\r)
\\=&\sum_{0<k<d}\f{(-1)^k}k\(\l(\f{d-1-k}3\r)+\l(\f{d-k}3\r)+\l(\f{d+1-k}3\r)\)
\\&+\f{(-1)^d}d\l(\f{d+1-d}3\r)
\\=&\f{(-1)^d}d.
\endalign$$
Therefore
$$\align&F(d-1)+F(d)+F(d+1)
\\=&\f12\(\bi{2p^a}{p^a-d}-\l(\f{p^a-d+1}3\r)-\l(\f{p^a-d}3\r)-\l(\f{p^a-d-1}3\r)\)
\\&+(-1)^{p-1}p^a\f{(-1)^d}d-p[p=3]\(\l(\f{d-1}3\r)+\l(\f d3\r)+\l(\f{d+1}3\r)\)
\\=&\f{p^a}{p^a-d}\bi{2p^a-1}{p^a-d-1}+(-1)^{p^a-d-1}\f{p^a}d.
\endalign$$

Clearly
$$\align&(-1)^{p^a-d-1}\bi{2p^a-1}{p^a-d-1}=\prod_{0<r<p^a-d}\l(1-\f{2p^a}r\r)
\\\eq&1-\sum_{0<r<p^a-d}\f{2p^a}r=2\(1-\sum_{0<r<p^a-d}\f{p^a}r\)-1
\\\eq&2(-1)^{p^a-d-1}\bi{p^a-1}{p^a-d-1}-1\ (\mo\ p^2).
\endalign$$
So, by the above we have
$$\align&F(d-1)+F(d)+F(d+1)
\\\eq&\f{p^a}{p^a-d}(-1)^{p^a-d-1}+(-1)^{p^a-d-1}\f{p^a}d
\\\eq&(-1)^{p^a-d-1}\f{p^a}d\cdot\f{p^a}{p^a-d}\eq0\ (\mo\ p^2)
\endalign$$
and
$$\align&F(d-1)+F(d)+F(d+1)
\\\eq&2\bi{p^a}{p^a-d}-(-1)^{p^a-d-1}\f{p^a}{p^a-d}+(-1)^{p^a-d-1}\f{p^a}d
\\\eq&2\bi{p^a}d-p^a\(\f{(-1)^{p^a-d}}d+(-1)^{p^a-1}\f{(-1)^{d}}{p^a-d}\)\ (\mo\ p^3).
\endalign$$
This concludes the proof. \qed

\medskip
\noindent{\it Proof of Theorem 1.1}.
(i) For $k=1,\ldots,p^a-1$, clearly
$$\bi{p^a}k=\f{p^a}k\prod_{0<j<k}\l(\f{p^a}j-1\r)\eq p^a\f{(-1)^{k-1}}k\ (\mo\ p^2).$$
Thus, with the help of Lemma 2.3, for $d\in\{p^a,p^a-1\}$ we have
$$\align F(d)=&(-1)^{p-1}p^a\sum_{k=1}^{p^a-1}\f{(-1)^k}k\l(\f{d-k}3\r)-p[p=3]\l(\f d3\r)
\\\eq&(-1)^{p-1}\sum_{k=1}^{p^a-1}\bi{p^a}k\l(\f{k-d}3\r)-p[p=3]\l(\f d3\r)
\\\eq&(-1)^{p-1}\l(1+(-1)^{p^a}-3[3\mid p^a]\r)\l(\f{-d-p^a}3\r)-p[p=3]\l(\f d3\r)
\\\eq&p[p=2]\l(\f{d+p^a}3\r)\eq p[p=2\ \&\ 3\nmid d+2^a]\ (\mo\ p^2).
\endalign$$
This proves (1.3) for $d=p^a,p^a-1$.

Assume that $0<d<p^a$. If
$$F(d)\eq p[p=2\ \&\ 3\nmid (-1)^a+d]\ (\mo\ p^2)$$
and
$$F(d+1)\eq p[p=2\ \&\ 3\nmid (-1)^a+d+1]\ (\mo\ p^2),$$
then by Lemma 3.1 we have
$$\align F(d-1)\eq&-F(d)-F(d+1)
\\\eq& -p[p=2]([3\nmid (-1)^a+d]+[3\nmid (-1)^a+d+1])
\\\eq&-p[p=2](2-[3\nmid (-1)^a+d-1])
\\\eq& p[p=2\ \&\ 3\nmid (-1)^a+d-1]\ (\mo\ p^2).
\endalign$$

By induction. we obtain from the above that (1.3) holds for any
$d=p^a,p^a-1,\ldots,0$.

(ii) By Lemma 2.2,
$$\align &F(0)+2F(1)
\\=&\f12\(\sum_{k=0}^{p^a-1}\bi{2k}k-\l(\f{p^a}3\r)\)+\sum_{k=0}^{p^a-1}\bi{2k}{k+1}-\l(\f{p^a-1}3\r)-2p[p=3]
\\\eq&p[p=2]+p^2[p=3]+\f12\l(1-\l(\f{p^a}3\r)\r)-\l(\f{p^a-1}3\r)-2p[p=3]
\\\eq&p[p=2]+\f{p^2}2[p=3]\ (\mo\ p^3).
\endalign$$
With the help of Lemma 2.4,
$$\align &F(p^a)+2F(p^a-1)
\\=&(-1)^{p-1}p^a\sum_{k=1}^{p^a-1}\f{(-1)^k}k\(\l(\f{p^a-k}3\r)+2\l(\f{p^a-1-k}3\r)\)
\\&-p[p=3]\(\l(\f{p^a}3\r)+2\l(\f{p^a-1}3\r)\)
\\=&(-1)^{p-1}p^a\sum_{k=1}^{p^a-1}\f{(-1)^k}k(1-3[3\mid p^a-k])+2p[p=3]
\\\eq&(-1)^p\l(p[p=2]+\f p2[p=3]\r)+2p[p=3]
\\\eq&p[p=2]+\f{p^2}2[p=3]\ (\mo\ p^3).
\endalign$$
Therefore
$$F(0)-F(p^a)\eq-2(F(1)-F(p^a-1))\ (\mo\ p^3).$$

When $p=2$ and $d\in\{1,\ldots,p^a-1\}$, we clearly have
$$\f{p^a}d-\f{p^a}{p^a-d}\eq-\f{p^a}d-\f{p^a}{p^a-d}=-\f{p^a}d\cdot\f{p^a}{p^a-d}\eq0\ (\mo\ p^2)$$
and hence
$$\f{p^a}d-\f{p^a}{p^a-d}\eq\f{p^a}{p^a-d}-\f{p^a}d\ (\mo\ p^3).$$
Thus, whether $p=2$ or not, by Lemma 3.1 we always have
$$F(d-1)+F(d)+F(d+1)\eq F(p^a-d-1)+F(p^a-d)+F(p^a-d+1)\ (\mo\ p^3)$$
whenever $d\in\{1,\ldots,p^a-1\}$.

Set $D(i)=F(i)-F(p^a-i)$ for $i=0,1,\ldots,p^a$. If $0\ls i\ls p^a-3$ then
$$D(i)+D(i+1)+D(i+2)\eq0\eq D(i+1)+D(i+2)+D(i+3)\ (\mo\ p^3)$$
and hence $D(i+3)\eq D(i)\ (\mo\ p^3)$. If $p=3$ then
$-D(0)=D(p^a)\eq D(0)\ (\mo\ p^3)$ and hence $D(0)\eq0\ (\mo\ p^3)$.
If $p^a\eq1\ (\mo\ 3)$ then
$$-D(0)=D(p^a)\eq D(1)\ (\mo\ p^3).$$
If $p^a\eq2\ (\mo\ 3)$ then
$$-D(0)=D(p^a)\eq D(2)\eq-D(0)-D(1)\ (\mo\ p^3).$$
As $D(0)+2D(1)\eq0\ (\mo\ p^3)$, we always have
$D(i)\eq0\ (\mo\ p^3)$ for $i=0,1$.
Therefore
$$D(i)\eq0\ (\mo\ p^3)\quad\t{for all}\ i=0,1,2,\ldots,p^a.$$
So (1.4) is valid and we are done. \qed

\heading{4. Proofs of Theorem 1.2 and Corollary 1.2}\endheading

\medskip
\noindent{\it Proof of Theorem 1.2}. In the case $m=n=0$, (1.11)
reduces to (1.6).

Below we assume $m>0$.
By the Chu-Vandermonde convolution identity (cf. [GKP, (5.22)]),
$$\bi{p^am+2k}{p^an+k+d}=\sum_{j\in\Z}\bi{p^am}{p^an-j}\bi{2k}{k+j+d}$$
for any $k\in\N$. Thus we have
$$\align&\sum_{k=0}^{p^a-1}\bi{p^am+2k}{p^an+k+d}-\bi{p^am}{p^an}\sum_{k=0}^{p^a-1}\bi{2k}{k+d}
\\=&\sum_{j>0}\bi{p^am}{p^an-j}\sum_{k=0}^{p^a-1}\bi{2k}{k+j+d}
+\sum_{j>0}\bi{p^am}{p^an+j}\sum_{k=0}^{p^a-1}\bi{2k}{k-j+d}
\\=&\sum_{j=1}^{p^a-1}\bi{p^am}{p^an-j}\sum_{k=0}^{p^a-1}\bi{2k}{k+j+d}
+\sum_{j=1}^{p^a-1}\bi{p^am}{p^an+j}\sum_{k=0}^{p^a-1}\bi{2k}{k+j-d}
\\&+\bi{p^am}{p^an+p^a}\sum_{k=0}^{p^a-1}\bi{2k}{k+p^a-d}+R_d,
\endalign$$
where
$$\align R_d=&\sum_{p^a<j<p^a+d}\bi{p^am}{p^an+j}\sum_{k=0}^{p^a-1}\bi{2k}{k+j-d}
\\=&\sum_{0<j<d}\bi{p^am}{p^a(n+1)+j}\sum_{k=0}^{p^a-1}\bi{2k}{k+j+p^a-d}.
\endalign$$
Note that $R_d=0$ if $d\in\{0,1\}$.

By Lemma 2.1, there are $p$-adic integers $u$ and $v$ such that
$$\bi{p^am}{p^an}=\bi mn\l(1+[p=2]p(m-n)n+p^2(m-n)u\r)$$
and
$$\align&\bi{p^am}{p^a(n+1)}=\bi m{n+1}\l(1+[p=2]p(n+1)(m-n-1)+p^2(n+1)v\r)
\\&\quad=\bi mn\l(\f{m-n}{n+1}+[p=2]p(m-n)(m-n-1)+p^2(m-n)v\r).
\endalign$$

Let $j\in\{1,\ldots,p^a-1\}$.
When $n\not=m$, we have
$$\align\bi{p^am}{p^an+j}=&\f{(p^am)!}{(p^an)!(p^am-p^an)!}
\times\f{\prod_{0\ls i<j}(p^am-p^an-i)}{(p^an+1)\cdots(p^an+j)}
\\=&\bi{p^am}{p^an}\f{p^a(m-n)}{p^an+j}\prod_{0<i<j}\f{p^a(m-n)-i}{p^an+i}
\endalign$$
and hence
$$\f{\bi{p^am}{p^an+j}}{(m-n)\bi{m}n}\eq\f{\bi{p^am}{p^an}}{\bi{m}n}
\cdot\f {p^a}j\prod_{0<i<j}\f{p^a-i}i\eq\bi {p^a}j\ (\mo\ p^2)$$
since
$$\f{p^a}j-\f{p^a}{p^an+j}=\f{p^a}j\cdot\f{p^an}{p^an+j}\eq0\ (\mo\ p^2)$$
and
$$\f{p^a(m-n)-i}{p^an+i}\eq\f{p^a-i}{p^an+i}\eq\f{p^a-i}i\ (\mo\ p)\quad\t{for}\ \ 0<i<p^a.$$
Similarly, if $n\not=0$ then
$$\f{\bi{p^am}{p^an-j}}{n\bi mn}=\f{\bi{p^am}{p^a(m-n)+j}}{(m-(m-n))\bi m{m-n}}\eq\bi {p^a}j\ (\mo\ p^2).$$
Also,
$$\align\f{n+1}{\bi mn}\bi{p^am}{p^a(n+1)+j}=&(m-n)\f{\bi{pm}{p^a(n+1)+j}}{\bi m{n+1}}
\\\eq&(m-n)(m-n-1)\bi {p^a}j\ (\mo\ p^{2+\ord_p(m-n)}).
\endalign$$

Combining the above, we find that
$$\aligned&\bi{m}n^{-1}\sum_{k=0}^{p^a-1}\bi{p^am+2k}{p^an+k+d}-(1+[p=2]pn(m-n))\sum_{k=0}^{p^a-1}\bi{2k}{k+d}
\\&-\f{m-n}{n+1}\sum_{k=0}^{p^a-1}\bi{2k}{k+p^a-d}-\f{R_d}{\bi mn}
\\\eq&\sum_{j=1}^{p^a-1}\bi {p^a}j\sum_{k=0}^{p^a-1}\(n\bi{2k}{k+j+d}+(m-n)\bi{2k}{k+j-d}\)
\ (\mo\ p^2)
\endaligned\tag4.1$$
and
 $$\f{n+1}{\bi mn}R_d\eq(m-n)(m-n-1)\sum_{0<j<d}\bi {p^a}j\sum_{k=0}^{p^a-1}\bi{2k}{k+j+p^a-d}\ (\mo\ p^2).\tag4.2$$

By (1.6),
$$\aligned\sum_{k=0}^{p^a-1}\bi{2k}{k+d}\eq&\l(\f{p^a-d}3\r)-p[p=3]\l(\f d3\r)+2p^aS_d\ (\mo\ p^2),
\\\sum_{k=0}^{p^a-1}\bi{2k}{k+p^a-d}\eq&\l(\f{d}3\r)-p[p=3]\l(\f {p^a-d}3\r)+2p^aS_{p^a-d}\ (\mo\ p^2).
\endaligned\tag4.3$$
Clearly $p^aS_d$ is congruent to
$$\sum_{0<k<d}\(\f{p^a}k\prod_{0<j<k}\f{p^a-j}j\)\l(\f{d-k}3\r)=\sum_{0<k<d}\bi{p^a}k\l(\f{d-k}3\r)$$
modulo $p^2$.
Observe that both $p^aS_{p^a-d}$ and
$$\sum_{j=1}^{p^a-1}\bi {p^a}j\sum_{k=0}^{p^a-1}\bi{2k}{k+j+d}$$
are congruent to
$$\sum_{0<j<p^a-d}\bi {p^a}j\l(\f{p^a-j-d}3\r)=\sum_{d<k<p^a}\bi {p^a}k\l(\f{k-d}3\r)$$
modulo $p^2$. Also,
$$\align&\sum_{d<k<p^a}\bi {p^a}k\l(\f{k-d}3\r)
=\sum_{k=1}^{p^a-1}\bi {p^a}k\l(\f{k-d}3\r)-\sum_{0<k\ls d}\bi {p^a}k\l(\f{k-d}3\r)
\\=&\l(1+(-1)^{p^a}-3[3\mid p^a]\r)\l(\f{-d-p^a}3\r)+\sum_{0<k\ls d}\bi {p^a}k\l(\f{d-k}3\r)
\endalign$$
with the help of Lemma 2.3. Thus,
$$\aligned&p^aS_{p^a-d}\eq\sum_{j=1}^{p^a-1}\bi {p^a}j\sum_{k=0}^{p^a-1}\bi{2k}{k+j+d}
\\\eq&p[p=2\ \&\ 3\nmid d+(-1)^a]+p[p=3]\l(\f d3\r)+p^aS_d\ (\mo\ p^2).
\endaligned\tag4.4$$
Note that
$$\align&\sum_{j=1}^{p^a-1}\bi {p^a}j\sum_{k=0}^{p^a-1}\bi{2k}{k+j-d}
\\=&\sum_{0<j<d}\bi {p^a}j\sum_{k=0}^{p^a-1}\bi{2k}{k+d-j}+\sum_{d\ls j<p^a}\bi {p^a}j\sum_{k=0}^{p^a-1}\bi{2k}{k+j-d}
\endalign$$
is congruent to
$$\align&\sum_{0<j<d}\bi {p^a}j\l(\f{p^a-d+j}3\r)+\sum_{d\ls j<p^a}\bi {p^a}j\l(\f{p^a-j+d}3\r)
\\=&\sum_{0<j<d}\bi {p^a}j\(\l(\f{p^a-d+j}3\r)-\l(\f{p^a-j+d}3\r)\)+\sum_{k=1}^{p^a-1}\bi {p^a}k\l(\f{k+d}3\r)
\\=&(1-p[p=3])\sum_{0<j<d}\bi {p^a}j\l(\f{d-j}3\r)+(1+(-1)^{p^a}-3[3\mid p^a])\l(\f{d-p^a}3\r)
\endalign$$
modulo $p^2$ in light of (1.6) and Lemma 2.3. Therefore
$$\aligned&\sum_{j=1}^{p^a-1}\bi{p^a}j\sum_{k=0}^{p^a-1}\bi{2k}{k+j-d}
\\\eq& p^aS_d+p[p=2\ \&\ 3\nmid d-(-1)^a]-p[p=3]\l(\f d3\r)\ (\mo\ p^2).
\endaligned\tag4.5$$
As
$$\align&\sum_{0<j<d}\bi{p^a}j\sum_{k=0}^{p^a-1}\bi{2k}{k+j+p^a-d}
\\\eq&\sum_{0<j<d}\bi{p^a}j\l(\f{p^a-(j+p^a-d)}3\r)\eq p^aS_d\ (\mo\ p^2),
\endalign$$
(4.2) yields that
$$\f{n+1}{\bi mn}R_d\eq(m-n)(m-n-1)p^aS_d\ (\mo\ p^2).\tag4.6$$

Combining (4.1) with (4.3)--(4.6), we finally obtain (1.11). Note
that (1.12) follows from (1.11) in the case $m=2n$. This ends the
proof. \qed

\medskip
\noindent{\it Proof of Corollary 1.2}. By (1.12) in the case
$d\in\{0,1\}$, we have
$$\align&\f1{C_n}\sum_{k=0}^{p^a-1}C_{p^an+k}=\f1{C_n}\sum_{k=0}^{p^a-1}\(\bi{2(p^an+k)}{p^an+k}
-\bi{2(p^an+k)}{p^an+k+1}\)
\\\eq&(n+1)\l(\f {p^a}3\r)-\l(n+(n+1)\l(\f{p^a-1}3\r)-[p=3]p(n+1)\r)
\\\eq&1-3(n+1)\l(\f{p^a-1}3\r)\ (\mo\ p^2).
\endalign$$
This proves (1.13). On the other hand, (1.12) in the case $d=0$
yields
$$\f1{C_n}\sum_{k=0}^{p^a-1}(p^an+k+1)C_{p^an+k}\eq(n+1)\l(\f{p^a}3\r)\ (\mo\ p^2).$$
So we have
$$\align\f1{C_n}\sum_{k=0}^{p^a-1}kC_{p^an+k}
\eq&(n+1)\l(\f{p^a}3\r)-\f{p^an+1}{C_n}\sum_{k=0}^{p^a-1}C_{p^an+k}
\\\eq&(n+1)\l(\f{p^a}3\r)-(p^an+1)\l(1-3(n+1)\l(\f{p^a-1}3\r)\r)
\\\eq&(3p^an+1)(n+1)\l(\f{p^a-1}3\r)-p^an-1
\\&+(n+1)\l(2\l(\f{p^a-1}3\r)+\l(\f{p^a}3\r)\r)\ (\mo\ p^2).
\endalign$$
Since
$$2\l(\f{p^a-1}3\r)+\l(\f{p^a}3\r)=1-p[p=3],$$
(1.14) follows at once. We are done. \qed

\enddocument